# Analysis and Design of Actuation-Sensing-Communication Interconnection Structures towards Secured/Resilient Closed-loop Systems

Sérgio Pequito † Farshad Khorrami ‡ Prashanth Krishnamurthy ‡ George J. Pappas †

**Abstract**

This paper considers the analysis and design of resilient/robust decentralized control systems. Specifically, we aim to assess how the pairing of sensors and actuators lead to architectures that are resilient to attacks/hacks for industrial control systems and other complex cyber-physical systems. We consider inherent structural properties such as internal fixed modes of a dynamical system depending on actuation, sensing, and interconnection/communication structure for linear discrete time-invariant dynamical systems. We introduce the notion of resilient fixed-modes free system that ensures the non-existence of fixed modes when the actuation-sensing-communication structure is compromised due to attacks by a malicious agent on actuators, sensors, or communication components and natural failures. Also, we provide a graph-theoretical characterization for the resilient structurally fixed modes that enables to capture the non-existence of resilient fixed modes for almost all possible systems' realizations. Additionally, we address the minimum actuation-sensing-communication co-design ensuring the non-existence of resiliently structurally fixed modes, which we show to be NP-hard. Notwithstanding, we identify conditions that are often satisfied in engineering settings and under which the co-design problem is solvable in polynomial-time complexity. Furthermore, we leverage the structural insights and properties to provide a convex optimization method to design the gain for a parametrized system and satisfying the sparsity of a given information pattern. Thus, exploring the interplay between structural and non-structural systems to ensure their resilience. Finally, the efficacy of the proposed approach is demonstrated on a power grid example.

## I. Introduction

Complex cyber-physical and industrial control systems (ICS) are prevalent in a wide range of application domains such as multi-agent networks, infrastructure systems such as the electric power grid, autonomous

This work was supported in part by the TerraSwarm Research Center, one of six centers supported by the STARnet phase of the Focus Center Research Program (FCRP) a Semiconductor Research Corporation program sponsored by MARCO and DARPA, and the NSF ECCS-1306128 grant.

†Department of Electrical and Systems Engineering, School of Engineering and Applied Science, University of Pennsylvania

‡Department of Electrical and Computer Engineering, Polytechnic Institute of New York University

vehicles, industrial machinery, process control, and manufacturing systems [1]–[4]. Modern ICS are complex interconnected combinations of hardware and software components including sensors, actuators, physical dynamical systems/processes, computational systems (e.g., microprocessors, programmable logic controllers – PLCs), and communication mechanisms. With the increasing complexity and connectivity of cyber-physical systems, the potential attack surface for cyber-attacks is also increasing [5]–[12].

There are, in general, various types of attacks (Figure 1) relevant in a cyber-physical control system application including:

- modification of the behavior (control parameters, logic, timing, or run-time process variables) of a computational node (e.g., a computer implementing a controller);
- modification of messages between computational nodes; and
- modification of sensor or actuator messages/signals.

To develop a truly cyber-resilient system, each computational node, sensor input, actuator output, and communication link must be considered as potentially vulnerable, and a resilient control methodology must be developed to accommodate malicious modifications, failures, and malfunctions of any subsets of the overall set of these components.

Current-day real-world systems often inherently are of a decentralized control structure due to constraints/requirements such as communication-related restrictions, requirements for fault tolerance and resiliency, security-related considerations, etc. Due to the decentralized nature of the sensing-actuation capabilities of ICS, there typically exist a multitude of decision-makers, each of which only has access to partial data. Some groundbreaking work in the understanding of necessary and sufficient conditions to ensure arbitrary spectrum placement of closed-loop systems constrained to specified information patterns can be found in [1], [13]–[15]. In recent years, research in decentralized control has seen several advances that contributed to a renewed interest in the field [16]–[22]. To address the cyber-resiliency needs, sensor-actuator-communication structural interconnection and decentralized control system design and implementation techniques are developed in this paper.

In this paper, our aim is to develop a methodology to co-design actuation-sensing-communication topologies to attain a given performance under disruptive scenarios such as natural failures or attacks by a malicious agent on actuators, sensors, or communication components of the system. The threat model considered here is the disabling of a subset of actuators, sensors, and sensor-actuator feedback links in a decentralized control system. This class of attacks can be modeled as zeroing out specific "attacked" elements of the closed-loop feedback gain matrix. Such an attack on a sensor-actuator feedback interconnection could, as illustrated in Figure 1, result from an attack on a sensor, a computational node, or an actuator.

To achieve resiliency against such attacks, the proposed design methodology seeks to guarantee certain closed-loop system properties regardless of attack/non-attack conditions. Specifically, we aim to determine the minimum sensor-actuator pairing structure and a corresponding control algorithm that together guarantee that the closed-loop system attains a given performance under disruptive scenarios. In particular, given a system, one aspect of the problem considered here is to determine the structure of sensor locations (i.e., state dependencies of sensor measurements), actuator locations (i.e., state variables directly controlled), and communication between sensors and actuators so as to provide certain resiliency properties of the resulting closed-loop system.

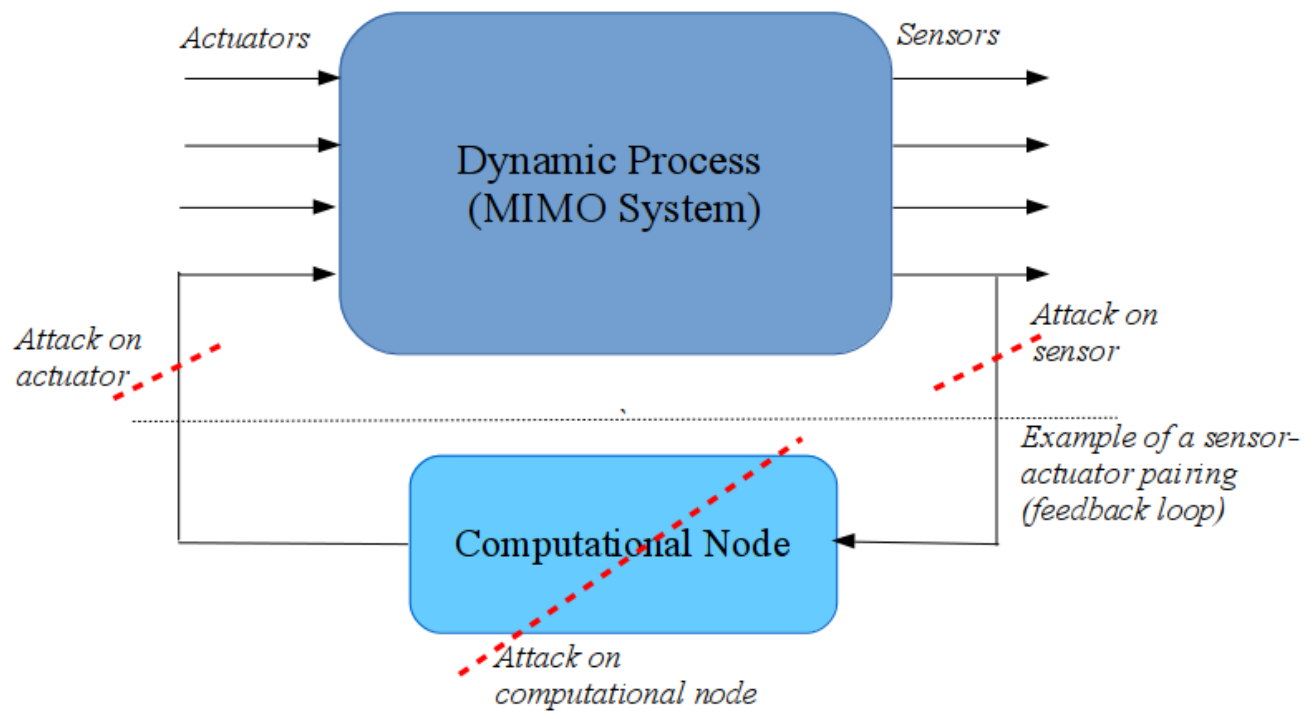

Fig. 1: Closed-loop feedback structure in an ICS and illustration of the various possible attack locations. The proposed methodology addresses a general multi-input-multi-output (MIMO) system with multiple decentralized feedback controllers; the figure shows one feedback link for simplicity.

Several necessary and sufficient conditions that characterize the structural controllability, as well as their verification, are known for linear time-invariant [23], and switching systems [24]. However, the design of actuation capabilities to ensure that these conditions hold has only been addressed in the last few years [25]–[31]. In particular, the problem of determining the minimum number of actuated variables ensuring structural controllability in LTI systems was addressed in [27], [29], and later extended to the case when cost constraints on the actuated state variables are imposed [28]. Further, when the actuation capabilities are known a priori, the problem of choosing the minimum subset of these to ensure structural controllability was shown to be NP-hard in [31]. The problem of determining the minimum number of actuated variables to ensure structural controllability was extended for discrete-time fractional dynamics [32] and linear time-invariant switching systems [33]. It was shown in [30] that determining the minimum actuation-sensing-communication required to ensure the existence of decentralized control laws in a generic structural system is NP-hard. However, none of these results address the design of secured/resilient decentralized closed-loop systems as we propose in the current paper. Yet, it is to be noted that the design of robust actuation/sensing configurations has been previously addressed and shown to be NP-hard [34]; more specifically, when the goal is to determine the minimum number of actuated state variables such that the system is structurally controllable under an actuator failure. The feasible minimal communication-pairings between sensors and actuators allowing decentralized control laws, i.e., essential information patterns, were studied in [35], and shown to be an NP-hard problem. In [36], the $p$-robustness of the power grid was analyzed with respect to transmission line failures, i.e., the number of failures that can occur while guaranteeing generic controllability of the associated dynamical system was determined. Hereafter, we leverage the insights gained in previous works and extend the previous literature to enable the actuation-sensing-communication co-design to cope with disruptive scenarios. To the best of authors' knowledge this is the first approach of this kind to the proposed co-design problem, which require the introduction of new concepts and methodologies, as we propose to in this paper.

The main contributions of this paper are fivefold: (*i*) introduction of the notion of resilient fixed-modes free system that ensures the non-existence of fixed modes when the actuation-sensing-communication structure is compromised; (*ii*) a graph-theoretical representation that ensures almost always the non-existence of resilient fixed modes; (*iii*) we show that the minimum actuation-sensing-communication co-design problem under disruptive scenarios, i.e., the minimum actuation-sensing-communication structure that ensures a system without resilient fixed modes almost always, is NP-hard; (*iv*) we propose an efficient (i.e., polynomial-time) solution to the minimum actuation-sensing-communication co-design problem under disruptive scenarios for subclasses of dynamical systems that are often found in practical engineering settings; and (*v*) determining, for a parametrized system, a stabilizing gain satisfying the sparsity of a given information pattern using convex optimization tools.

The remaining of this paper is organized as follows. The formal problem statement is provided in Section II. Section III reviews some concepts and introduces results in structural systems theory and computational complexity. The main technical results are presented in Section IV, followed by an illustrative example in Section V. Concluding remarks are presented in Section VI.

## II. Problem Formulation

Consider a (possibly) large-scale system modeled by

$$x[k+1] = Ax[k] + Bu[k], \quad y[k] = Cx[k], \tag{1}$$

where $x \in \mathbb{R}^n$ is the state of the system, $u \in \mathbb{R}^p$ is the input vector collecting the control signals and $y \in \mathbb{R}^m$ are the output vector collecting the measurements. To reduce the computational requirements to determine the input signal at the controller, one may consider a static output feedback (SOF) strategy provided as follows:

$$u[k] = Ky[k] = KCx[k], \tag{2}$$

hence, the *closed-loop control system* (1) using SOF (2) is described as follows:

$$x[k+1] = (A + BKC)x[k]. \tag{3}$$

For notational convenience, in the remainder of the paper, we refer to (3) by the tuple $(A, B, C, K)$, and we identify the system in (1) with the tuple $(A, B, C)$. However, the data collected by the sensors is not always available to all actuators, which implies that the matrix $K$ may have some zeros, i.e., it is not a full matrix. To formalize these notions, which are the core of decentralized control, we need to account for the *information patterns*, i.e., the sparsity of the *gain* matrix $K$, given by a binary matrix $\bar{K} \in \{0, 1\}^{p\times m}$, where $\bar{K}_{i,j} = 1$ if the measurements from sensor $j$ are available by actuator $i$, and zero otherwise.

*Remark 1:* For simplicity in illustrating the primary concepts introduced in this paper, the closed-loop system is formulated above as an LTI dynamic system and a static output feedback controller. Observe that the basic proposed design concepts can be extended to account for dynamic controllers that keeps track of the history or includes an observer to estimate unmeasured parts of the state, instead of simply constructing a feedback as a static gain applied to measured outputs [1], [17], [18]. ◇

Now, suppose that a natural failure or an attack to the actuator-sensor-communication structure occurs. Then, we need the notion of *robust information patterns*, i.e., the actuation-sensing-communication topology, associated with the existence of decentralized control strategies under disruptive scenarios. Subsequently, we introduce the notion of *resilient fixed modes* that are a natural extension of fixed modes [13] (i.e., the eigenvalues of the closed-loop system that cannot be changed by considering different gains satisfying a prescribed information pattern), when accounting for communication failures, or, equivalently, actuators and sensor failures as explained before (see Figure 1).

*Definition 1:* Consider system (1)-(2) with $p$ actuators and $m$ sensors, $\mathcal{K}$ the set of possible communication links from the sensors to the actuators, and $\{\Gamma_i\}_{i\in\Delta}$ the collection of possible communication failures, with $\Gamma_i \subset \mathcal{K}$ for all $i \in \Delta$. The set

$$\Sigma(A, B, C; \mathcal{K}, \{\Gamma_i\}_{i\in\Delta}) = \bigcup_{i\in\Delta} \bigcap_{K\in[\bar{K}(\mathcal{K})]^{\Gamma_i}_{p,m}} \sigma(A + BKC)$$

is defined to be the set of *resilient* fixed modes under a collection of feedback failures $\{\Gamma_i\}_{i\in\Delta}$ of the closed-loop system (3) with respect to the information pattern $\bar{K}(\mathcal{K})$, where $\bar{K}_{i,j} = 1$ if $(j, i) \in \mathcal{K}$ and zero otherwise, and $[\bar{K}(\mathcal{K})]^{\Gamma_i}_{p,m}$ is the set of all possible $p \times m$ constant output feedback matrices given by $[\bar{K}(\mathcal{K})]^{\Gamma_i}_{p,m} = \{K \in \mathbb{R}^{p\times m} : \ K_{lk} = 0 \text{ if } (k, l) \notin \mathcal{K} \text{ or } (k, l) \in \Gamma_i, \ 0 \le l \le p, \, 0 \le k \le m\}$. ◇

Notice that by ensuring that the system does not have any resilient fixed modes, it is guaranteed that under any resulting information pattern, all the closed-loop system eigenvalues can be arbitrarily placed by choosing the feedback gain appropriately.

Additionally, in many practical scenarios, it is often the case that the exact values of the non-zero parameters of the plant matrices are unknown. To address this problem, in this paper, we adopt the framework of structural systems [23]. To this end, let $\bar{A} \in \{0,1\}^{n\times n}$, $\bar{B} \in \{0,1\}^{n\times p}$ and $\bar{C} \in \{0,1\}^{n\times m}$ be the binary matrices that represent the structural patterns (location of fixed zeros and unknown parameters) of $A, B$, and $C$, respectively. We then focus on properties of systems, where the plant matrices have these sparsity patterns $(\bar{A}, \bar{B}, \bar{C})$ which we refer to as a *structural system*. In the context of structural systems, controllability and observability are referred to as structural controllability and structural observability, respectively. The structural systems counterpart of the fixed modes concept is referred to as *structural fixed modes*, which, essentially, are the fixed modes attributed to the structural pattern of a system, as opposed to those that originate from the canceling of numerical parameters [37], [38]. Therefore, almost all systems in the sparsity class $(\bar{A}, \bar{B}, \bar{K}, \bar{C})$ have no fixed modes, and, hence, allow pole placement arbitrarily close to any pre-specified (symmetrical about the real axis) set of eigenvalues by a static output feedback with the sparsity of $\bar{K}$ [37]. Subsequently, the resilient fixed-modes can also be understood from a structural systems perspective, and formally defined as follows.

*Definition 2:* Consider system (1)-(2) with $p$ actuators and $m$ sensors, and denote its structural pattern by $(\bar{A}, \bar{B}, \bar{C}; \mathcal{K}, \{\Gamma_i\}_{i\in\Delta})$, where $\mathcal{K}$ is the set of possible communication links from the sensors to the actuators, and $\{\Gamma_i\}_{i\in\Delta}$ the collection of possible communication failures, with $\Gamma_i \subset \mathcal{K}$ for all $i \in \Delta$. Then, $(\bar{A}, \bar{B}, \bar{C}; \mathcal{K}, \{\Gamma_i\}_{i\in\Delta})$ is said to have structurally resilient fixed modes w.r.t. an information pattern $\bar{K}(\mathcal{K})$, where $\bar{K}_{i,j} = 1$ if $(j,i) \in \mathcal{K}$ and zero otherwise, if for almost all $A \in [\bar{A}]_{n\times n}$, $B \in [\bar{B}]_{n\times p}$, $C \in [\bar{C}]_{m\times n}$, we have $\Sigma(A, B, C; \mathcal{K}, \{\Gamma_i\}_{i\in\Delta}) \neq \emptyset$. ◇

Hereafter, we focus on the *analysis* and *design* of secured/resilient closed-loop control systems. Therefore, the first problem addresses the analysis, and can be formulated as follows:

$\mathcal{P}_1$ : Provide necessary and sufficient conditions to ensure that a system has no structurally resilient fixed modes w.r.t. an information pattern and a collection of possible communication failures. ○

Then, with these conditions, we aim to design *secured/resilient* large-scale closed-loop control systems. More precisely, suppose that only the structure of the autonomous system is known, i.e., $\bar{A}$. Then, the goal is to determine the location of *dedicated* actuators and sensors (i.e., that only actuate and measure a single state variable, respectively). This scenario is commonly found in practice since state variables have specific physical properties associated with them [2], [3]. For instance, in power systems one often controls/measure the frequency of different generators and aggregate loads [39]. In addition, the communication between sensors and actuators, such that the system has no structurally resilient fixed modes w.r.t. the information pattern induced by the chosen communication links and any set of possible actuator-sensor-communication failures containing at most $k$ elements. In other words, given $\bar{A}$, we want to determine $(\bar{B}(\mathcal{I}_B), \bar{C}(\mathcal{I}_C), \bar{K}(\mathcal{K}))$, where $\bar{B}$ (respectively, $\bar{C}$) contain columns (respectively, rows) with at most one nonzero entry in each column (respectively, row) indexed by the elements in the collection $\mathcal{I}_B \equiv \{\iota_i\}_{i\in\mathcal{I}}$ with $\iota_i \in \{1, \cdots, n\}$ (respectively, $\mathcal{I}_C \equiv \{\iota_j\}_{j\in\mathcal{J}}$), and $K_{i,j} = 1$ if $(j,i) \in \mathcal{K}$, and zero otherwise. Notice that, in particular, $\mathcal{I}_B$ and $\mathcal{I}_C$ may contain several times the same index, which corresponds to the case where the same state variable is actuated and measured by different actuators and sensors, respectively. On the contrary, $\mathcal{K}$ does not contain twice the same pair of indices since these correspond uniquely to a sensor-actuator pair, and multiple simultaneous communications between the sensor and actuator in the pair are not allowed. Therefore, we have the following reslient actuation-sensing-communication co-design problem:

$\mathcal{P}_2$ : Given $\bar{A}$, and a maximum of $k$ actuator-sensor-communication failures, we want to determine $(\bar{B}(\mathcal{I}_B^*), \bar{C}(\mathcal{I}_C^*), \bar{K}(\mathcal{K}^*))$, where $(\mathcal{I}_B^*, \mathcal{I}_C^*, \mathcal{K}^*)$ is a solution to the following optimization problem:

$$\begin{aligned} \min_{\mathcal{I}_B, \mathcal{I}_C, \mathcal{K}} \quad & |\mathcal{I}_B| + |\mathcal{I}_C| + |\mathcal{K}| \\ \text{s. t.} \quad & (\bar{A}, \bar{B}(\mathcal{I}_B \setminus \mathcal{J}_B), \bar{C}(\mathcal{I}_C \setminus \mathcal{J}_C); \mathcal{K}, \mathcal{J}) \text{ has no struct.} \\ & \text{resilient fixed modes, for all } \mathcal{J} \subset \mathcal{K}, \mathcal{J}_B \subset \mathcal{I}_B, \\ & \mathcal{J}_C \subset \mathcal{I}_C, \text{ and } |\mathcal{J}| + |\mathcal{J}_B| + |\mathcal{J}_C| \le k, \end{aligned} \tag{4}$$

where $|\mathcal{I}_B|$ and $|\mathcal{I}_C|$ are the number of dedicated actuators and sensors, respectively, i.e., the number of elements in the collections $\mathcal{I}_B$ and $\mathcal{I}_C$, and $|\mathcal{K}|$ are the number of communication links from the sensors to the actuators to obtain the decentralized closed-loop control system. ◦

Unfortunately, the co-design problem $\mathcal{P}_2$ is NP-hard, as we show in the Theorem 2. Therefore, in this paper we identify a set of conditions that are commonly found in engineering applications (see Remark 3), and under which we provide a computationally tractable procedure.

Lastly, in many ICS applications, one is only required to ensure stabilizability of the closed-loop system, which implies that only the open-loop unstable eigenvalues need to be changed. In other words, it is not required to ensure placing of *all* the eigenvalues, which is known to be NP-hard [40]. Therefore, the last problem addressed in this paper can be stated as follows:

$\mathcal{P}_3$ : Given a realization $(A, B, C)$ of $(\bar{A}, \bar{B}, \bar{C})$ and $\bar{K}$ from the design in $\mathcal{P}_2$, determine a gain $K \in [\bar{K}]$ such that $(A, B, C, K)$ is asymptotically stable. ◦

Subsequently, we formulate an optimization problem whose feasibility allows to determine if the generic design obtained in $\mathcal{P}_2$ holds for a specific parameterization. Furthermore, we propose a convex optimization formulation whose feasibility enables us to retrieve a solution to $\mathcal{P}_3$.

## III. Preliminaries and terminology

In order to perform structural analysis efficiently, it is customary to associate to (1) a directed graph (digraph) $\mathcal{D} = (\mathcal{V}, \mathcal{E})$, in which $\mathcal{V}$ denotes the set of *vertices* and $\mathcal{E} \subseteq \mathcal{V} \times \mathcal{V}$ the set of *edges*, where $(v_j, v_i)$ represents an edge from the vertex $v_j$ to vertex $v_i$. Consider the structural systems $(\bar{A}, \bar{B}, \bar{C})$, then we can denote by $\mathcal{X} = \{x_1, \ldots, x_n\}$, $\mathcal{U} = \{u_1, \ldots, u_p\}$ and $\mathcal{Y} = \{y_1, \ldots, y_m\}$ the sets of state, input and output vertices respectively, and $\mathcal{E}_{\mathcal{X},\mathcal{X}} = \{(x_i, x_j) : \ \bar{A}_{ji} \neq 0\}$, $\mathcal{E}_{\mathcal{U},\mathcal{X}} = \{(u_j, x_i) : \ \bar{B}_{ij} \neq 0\}$, and $\mathcal{E}_{\mathcal{X},\mathcal{Y}} = \{(x_i, y_j) : \bar{C}_{ji} \neq 0\}$ the edges between the sets in subscript. In addition, given an *information pattern* $\bar{K} \in \{0,1\}^{p \times m}$, describing output feedback in the inputs, we also have $\mathcal{E}_{\mathcal{Y},\mathcal{U}} = \{(y_j, u_i) : \ \bar{K}_{ij} \neq 0\}$. In particular, we have the *state* digraph $\mathcal{D}(\bar{A}) = (\mathcal{X}, \mathcal{E}_{\mathcal{X},\mathcal{X}})$, the *system* digraph $\mathcal{D}(\bar{A}, \bar{B}, \bar{C}) = (\mathcal{X} \cup \mathcal{U} \cup \mathcal{Y}, \mathcal{E}_{\mathcal{X},\mathcal{X}} \cup \mathcal{E}_{\mathcal{U},\mathcal{X}} \cup \mathcal{E}_{\mathcal{X},\mathcal{Y}})$, and the *closed-loop system* digraph $\mathcal{D}(\bar{A}, \bar{B}, \bar{C}, \bar{K}) = (\mathcal{X} \cup \mathcal{U} \cup \mathcal{Y}, \mathcal{E}_{\mathcal{X},\mathcal{X}} \cup \mathcal{E}_{\mathcal{U},\mathcal{X}} \cup \mathcal{E}_{\mathcal{X},\mathcal{Y}} \cup \mathcal{E}_{\mathcal{Y},\mathcal{U}})$.

In what follows, we consider the well known concepts of *directed path*, *elementary path*, *cycle*, *strongly connected component* (SCC) and *directed acyclic graph* (DAG) [27], [41]. In addition, we can classify the SCCs in a DAG with respect to the existence of incoming and/or outgoing edges as follows: an SCC is *non-top linked* if it has no incoming edges from another SCC, and *non-bottom linked* if it has no outgoing edges to another SCC.

For any digraph $\mathcal{D} = (\mathcal{V}, \mathcal{E})$ and any two vertex sets $\mathcal{S}_1, \mathcal{S}_2 \subset \mathcal{V}$, we define the *bipartite graph* $\mathcal{B}(\mathcal{S}_1, \mathcal{S}_2, \mathcal{E}_{\mathcal{S}_1,\mathcal{S}_2})$ whose vertex set is given by $\mathcal{S}_1 \cup \mathcal{S}_2$ (where $\mathcal{S}_1$ and $\mathcal{S}_2$ are assumed to be disjoint) and the edge set $\mathcal{E}_{\mathcal{S}_1,\mathcal{S}_2} = \mathcal{E} \cap (\mathcal{S}_1 \times \mathcal{S}_2)$. In addition, we will make heavy use of the *state bipartite graph* $\mathcal{B}(\bar{A}) \equiv \mathcal{B}(\mathcal{X}, \mathcal{X}, \mathcal{E}_{\mathcal{X},\mathcal{X}})$, which is the bipartite graph *associated* with the state digraph $\mathcal{D}(\bar{A}) = (\mathcal{X}, \mathcal{E}_{\mathcal{X},\mathcal{X}})$.

Given $\mathcal{B}(\mathcal{S}_1, \mathcal{S}_2, \mathcal{E}_{\mathcal{S}_1,\mathcal{S}_2})$, a matching $M$ corresponds to a subset of edges in $\mathcal{E}_{\mathcal{S}_1,\mathcal{S}_2}$ so that no two edges have a vertex in common. A maximum matching $M^*$ is a matching $M$ that has the largest number of edges among all possible matchings. We call the vertices in $\mathcal{S}_1$ and $\mathcal{S}_2$ belonging to an edge in $M^*$, the *matched vertices* with respect to (w.r.t.) $M^*$, and *unmatched vertices* otherwise. If the set of unmatched vertices associated with a maximum matching are empty, then we have a *perfect matching*. For ease of referencing, in the remaining of the paper, the term *right-unmatched vertices* associated with the matching $M$ of $\mathcal{B}(\mathcal{S}_1, \mathcal{S}_2, \mathcal{E}_{\mathcal{S}_1,\mathcal{S}_2})$ (not necessarily maximum), will refer to those vertices in $\mathcal{S}_2$ that do not belong to a matching edge in $M^*$, dually a vertex from $\mathcal{S}_1$ that does not belong to an edge in $M^*$ is called a *left-unmatched vertex*, see, for instance, [27] for an illustrative example of the concepts.

Finally, we need to introduce the following definition to create the pairing between outputs and inputs.

*Definition 3 (Sequential-pairing):* Consider two sets of indices $\mathcal{I} = \{i_1, \ldots, i_n\}$ and $\mathcal{J} = \{j_1, \ldots, j_m\}$, and a maximum matching $M$ of the bipartite graph $\mathcal{B}(\mathcal{J}, \mathcal{I}, \mathcal{E}_{\mathcal{J},\mathcal{I}})$, where $\mathcal{E}_{\mathcal{J},\mathcal{I}} \subseteq \mathcal{J} \times \mathcal{I}$. We denote by

$|\mathcal{I},\mathcal{J}\rangle_M$ a sequential-pairing induced by $M$, defined as follows:

$$|\mathcal{I},\mathcal{J}\rangle_M = \left(\bigcup_{l=2,\ldots,k}\{(i_l,j_{l-1})\}\right)\cup\{(i_1,j_k)\},$$

where $(i_l,j_l)\in M$, for $l=1,\ldots,k$. $\diamond$

*Remark 2:* The sequential-pairing consists in the collection of edges such that $M\cup|\mathcal{I},\mathcal{J}\rangle_M$ forms a cycle. $\diamond$

## IV. Resilient actuation-sensing-communication co-design

In this section, we present the main results of this paper. First, in Theorem 1, we introduce the graph-theoretical characterization of resilient structurally fixed modes, which provide us with a characterization for analysis of resilient systems as solicited in $\mathcal{P}_1$. Then, in Theorem 2, we show that the resilient co-design problem presented in $\mathcal{P}_2$ is NP-hard. Notwithstanding, we characterize the solutions to $\mathcal{P}_2$, under assumptions that are often found in engineering settings, and that can be efficiently obtained (Theorem 4). Towards this goal, we explore different scenarios that enable us to gain some insights on how the structure constraints the existence of resilient structurally fixed modes. In particular, we show how that communication infrastructure can be established to ensure the resilience of the system when only communication failures/attacks are performed (see Theorem 3). Furthermore, in Section IV-A, we discuss how these assumptions can be waived in the case where we seek system's stabilization. Subsequently, in Section IV-B, we provide a feasibility criteria for the stabilization of a particular plant without unstable resilient fixed modes (Theorem 5), and an iterative optimization scheme is provided in Algorithm 1 to obtain a specific gain, which feasibility is detailed in Theorem 6; thus, providing a bridge between the structural design the existence of a specific gain when a specific parameterization of the system plant is considered.

We start by addressing $\mathcal{P}_1$, i.e., we provide necessary and sufficient conditions required to attain a closed-loop system that is free of resilient structurally fixed modes.

*Theorem 1:* The structural system $(\bar{A},\bar{B},\bar{C};\mathcal{K},\{\Gamma_i\}_{i\in\Delta})$ associated with (3) has no resilient structurally fixed modes if and only if both of the following conditions hold for all $i\in\Delta$:

$(a)$ in $\mathcal{D}(\bar{A},\bar{B},\bar{C};\mathcal{K},\Gamma_i)=(\mathcal{X}\cup\mathcal{U}\cup\mathcal{Y},\mathcal{E}_{\mathcal{X},\mathcal{X}}\cup\mathcal{E}_{\mathcal{X},\mathcal{Y}}\cup\mathcal{E}_{\mathcal{U},\mathcal{X}}\cup(\mathcal{E}_{\mathcal{Y},\mathcal{U}}\setminus\mathcal{E}_{\Gamma_i}))$, with $\mathcal{E}_{\Gamma_i}=\{(y_i,u_j):(j,i)\in\Gamma_i\}$, each state vertex $x\in\mathcal{X}$ is contained in an SCC which includes an edge of $\mathcal{E}_{\mathcal{Y},\mathcal{U}}\setminus\mathcal{E}_{\Gamma_i}$;

$(b)$ there exists a disjoint union of cycles $\mathcal{C}_k=(\mathcal{V}_k,\mathcal{E}_k)$ (subgraph of $\mathcal{D}(\bar{A},\bar{B},\bar{C};\ \mathcal{K},\Gamma_i)$) with $k\in\mathbb{N}$ such that $\mathcal{X}\subset\bigcup_{j=1}^{k}\mathcal{V}_j$. $\diamond$

These conditions are key in ensuring feasibility of the solutions attained in the resilient actuation-sensing-communication co-design problem. Nonetheless, these are also intrinsic to the computational complexity associated with determining a solution to $\mathcal{P}_2$, as we formally state in the next result.

*Theorem 2:* The resilient co-design problem $\mathcal{P}_2$ is NP-hard. ◇

Although $\mathcal{P}_2$ is NP-hard, it does not mean that all instances of the problem are equally computationally challenging. Therefore, in what follows we seek to determine conditions under which the problem is computationally tractable. Specifically, we propose to address $\mathcal{P}_2$ under the following assumption:

**A1** The structural dynamic matrix $\bar{A}$ is such that the state bipartite graph $\mathcal{B}(\bar{A}) = \mathcal{B}(\mathcal{X}, \mathcal{X}, \mathcal{E}_{\mathcal{X},\mathcal{X}})$, associated with $\bar{A}$, has a perfect matching. ◇

In other words, assumption **A1** can be interpreted as $\mathcal{D}(\bar{A})$ being spanned by a disjoint union of cycles (see [27] for details). In particular, this assumption is of practical interest in engineering applications as emphasized in the next remark.

*Remark 3:* We notice that this assumption covers a practically relevant class of systems since it subsumes the following scenarios: (i) almost all dynamical systems modeled as random networks of the Erdős-Rényi type [42]; (ii) discretized dynamical systems using Euler discretization for almost all discretization steps; (iii) several dynamical systems as presented in [2], [3]. ◇

Now, we consider the minimum actuation-sensing capabilities, and describe how to determine the sparsest information pattern that ensures the closed-loop system has no structurally resilient fixed modes. Towards this goal, we first consider the following procedure.

**Procedure 1**: Consider $\mathcal{I}_B^*$ and $\mathcal{I}_C^*$ to be the minimum collections of state variables' indices associated with dedicated actuators and sensors, respectively, such that the pair $(\bar{A}, \bar{B}(\mathcal{I}_B^*))$ and $(\bar{A}, \bar{C}(\mathcal{I}_C^*))$ is structurally controllable and observable, respectively. Further, without loss of generality, assume that $|\mathcal{I}_B^*| \leq |\mathcal{I}_C^*|$, with $\mathcal{I}_B = \{1, \ldots, |\mathcal{I}_B^*|\}$ and $\mathcal{I}_C = \{1, \ldots, |\mathcal{I}_C^*|\}$ the sets of indices numbering the dedicated inputs and outputs, respectively. Further, define the bipartite graph $\mathcal{B} = \mathcal{B}(\mathcal{I}_B, \mathcal{I}_C, \mathcal{E}_{\mathcal{I}_B,\mathcal{I}_C})$ having an edge $(i,j) \in \mathcal{E}_{\mathcal{I}_B,\mathcal{I}_C}$ if there is a path from dedicated actuator $i \in \mathcal{I}_B$ to dedicated sensor $j \in \mathcal{I}_C$ (in which case we say that an actuator *reaches* a sensor, or a sensor is *reached* by an actuator). Therefore, a maximum matching $M^*$ of $\mathcal{B}$ describes the maximum number of actuators that reach an equal number of sensors. Subsequently, by considering an information pattern $\bar{K}(\mathcal{K})$, where $\mathcal{K}$ is given by the sequential-pairing $|\mathcal{I}_B, \mathcal{I}_C\rangle_{M^*}$, it follows that $\mathcal{D}(\bar{A}, \bar{B}(\mathcal{I}_B^*), \bar{C}(\mathcal{I}_C^*), \bar{K}(\mathcal{K}))$ has a single non-top linked SCC that contains all non-top and non-bottom SCCs of $\mathcal{D}(\bar{A})$, with variables controlled and measured by the actuators and sensors, respectively, labelled by the indexes that belong to edges in $M^*$. Furthermore, it also contains all SCCs that are simultaneously reachable and reach these non-top and non-bottom SCCs. Intuitively, several SCCs are contracted into a single one with the auxiliary of feedback links, which ensure that the newly created non-top linked SCC satisfies Theorem 1-(a).

Nonetheless, $\mathcal{D}(\bar{A}, \bar{B}(\mathcal{I}_B^*), \bar{C}(\mathcal{I}_C^*), \bar{K}(\mathcal{K}))$ still has $|\mathcal{I}_C^*| - |\mathcal{I}_B^*|$ non-bottom linked SCCs, as well as (possibly) additional SCCs that reach these. Therefore, a total of $|\mathcal{I}_C^*| - |\mathcal{I}_B^*|$ additional feedback links, corresponding to the indices of $\mathcal{I}_C$ not present in the edges of $M^*$, need to be considered. Let $\mathcal{K}'$ correspond to feedback links starting in the dedicated sensors labelled by the former indices, and ending in only one of the dedicated actuators indexed by $\mathcal{I}_B$. Then, $\mathcal{D}(\bar{A}, \bar{B}(\mathcal{I}_B^*), \bar{C}(\mathcal{I}_C^*), \bar{K}(\mathcal{K} \cup \mathcal{K}'))$, consists in a single SCC with feedback links on it, which implies that the closed-loop digraph readily satisfies Theorem 1-(a). Thus, under the assumption **A1**, $\bar{K}(\mathcal{K}^*)$ with $\mathcal{K}^* = \mathcal{K} \cup \mathcal{K}'$ is a sparsest information pattern ensuring a closed-loop system without structurally fixed modes. $\circ$

Subsequently, we obtain the following result.

*Lemma 1:* Given $(\bar{A}, \bar{B}(\mathcal{I}_B^*), \bar{C}(\mathcal{I}_C^*))$, where $\mathcal{I}_B^*$ and $\mathcal{I}_C^*$ are the minimum collections of state variables ensuring structural controllability and observability, then Procedure 1 originates a sparsest information pattern $\bar{K}(\mathcal{K}^*)$ such that $(\bar{A}, \bar{B}(\mathcal{I}_B^*), \bar{C}(\mathcal{I}_C^*), \bar{K}(\mathcal{K}^*))$ has no structurally fixed modes, under assumption **A1**. In addition, Procedure 1 can be implemented in polynomial time, more precisely, $\mathcal{O}(\max\{|\mathcal{V}|+|\mathcal{E}|, (|\mathcal{I}_B^*|+|\mathcal{I}_C^*|)^{\frac{5}{2}}\})$, where $\mathcal{V}$ and $\mathcal{E}$ are the number of vertices and edges of $\mathcal{D}(\bar{A}, \bar{B}(\mathcal{I}_B^*), \bar{C}(\mathcal{I}_C^*))$, respectively. $\diamond$

In fact, we can consider the scenario where only failures/attacks at the communication infrastructure level occur. In such scenario, we need to provide robust information patterns, i.e., that can account for an arbitrary number $k$ of communication failure while ensuring a resilient structurally fixed modes free system, which can be obtained by extending Lemma 1 as described in the next result.

*Theorem 3:* Given $(\bar{A}, \bar{B}(\mathcal{I}_B^*), \bar{C}(\mathcal{I}_C^*))$, where $\mathcal{I}_B^*$ and $\mathcal{I}_C^*$ are the minimum collections of state variables' indices ensuring structural controllability and observability respectively, and $k$ the maximum number of communication failures, then $\bar{K}(\mathcal{K}^*)$ such that $\mathcal{K}^* = \mathcal{K}^{(1)} \cup \ldots \cup \mathcal{K}^{(k+1)}$, where $\mathcal{K}^{(i)}$ is determined using Procedure 1 and $\mathcal{K}^{(i)} \cap \mathcal{K}^{(j)} = \emptyset$ for $i \neq j$, is a sparsest information pattern that ensures $(\bar{A}, \bar{B}(\mathcal{I}_B^*), \bar{C}(\mathcal{I}_C^*), \bar{K}(\mathcal{K}^*))$ to have no structurally fixed modes w.r.t. $\bar{K}(\mathcal{K}^*)$ under assumption **A1**, and any set of failures indexed by $\mathcal{J} \subset \mathcal{K}^*$, with $|\mathcal{J}| \leq k$. In addition, these can be determined in polynomial time, more precisely, $\mathcal{O}(k \max\{|\mathcal{V}| + |\mathcal{E}|, (|\mathcal{I}_B^*| + |\mathcal{I}_C^*|)^{\frac{5}{2}}\})$, where $|\mathcal{V}|$ and $|\mathcal{E}|$ are the number of vertices and edges of $\mathcal{D}(\bar{A}, \bar{B}(\mathcal{I}_B^*), \bar{C}(\mathcal{I}_C^*))$, respectively. $\diamond$

*Remark 4:* To ensure that two sparsest information patterns using Procedure 1 are such that $\mathcal{K}^{(i)} \cap \mathcal{K}^{(j)} = \emptyset$ for $i \neq j$, one just needs to ensure that the edges in the maximum matching $M^{(i)}$ found to determine $\mathcal{K}^{(i)}$ are not considered when computing $\mathcal{K}^{(j)}$. Therefore, when computing $\mathcal{K}^{(j)}$, instead of the set $\mathcal{E}_{\mathcal{I}_B, \mathcal{I}_C}$ of the bipartite graph we just need to consider $\mathcal{E}_{\mathcal{I}_B, \mathcal{I}_C} \setminus \left( \bigcup_{l=1}^{j-1} M^{(l)} \right)$. Hence, the procedure does not increase the computational complexity. $\diamond$

Therefore, it is worth to notice that the sparsest communication required under $k$ possible communication failures does not consist in the direct union of $k+1$ (arbitrary) sparsest information patterns. Besides, the procedure to determine disjoint sparsest information patterns establishes an upper bound on the possible number of failures/attacks on the communication, since it might not be possible to determine additional sparsest information patterns that are disjoint.

Next, to characterize the solutions to $\mathcal{P}_2$, we need the following lemma, which accounts for the minimum collections of dedicated actuators that ensure structural controllability when $k$ failures might occur.

*Lemma 2:* Under assumption **A1**, the pair $(\bar{A}, \bar{B}(\mathcal{I}_B \setminus \mathcal{J}_B))$ is structurally controllable, with $\mathcal{J}_B \subset \mathcal{I}_B$ and $|\mathcal{J}_B| \leq k$ if and only if for each non-top linked SCC of $\mathcal{D}(\bar{A})$ there exists at least $k+1$ indices in $\mathcal{I}_B$ corresponding to its state variables. In particular, the minimum number of dedicated actuators required equals $(k+1)t$, where $t$ is the number of non-top linked SCCs, i.e., exactly $k+1$ actuators controlling state variables in each non-top linked SCC. ◇

Similarly, by duality between controllability and observability in LTI systems [43], we obtain a similar result regarding the minimum number of sensors required to attain observability under at most $k$ failures, which we refer to the *dual* of Lemma 2.

Therefore, a solution to $\mathcal{P}_2$ can be characterized as follows.

*Theorem 4:* Let $\mathcal{I}_B$ and $\mathcal{I}_C$ be the indices associated minimum collection of dedicated actuators and sensors $\bar{B}(\mathcal{I}_B^*)$ and $\bar{C}(\mathcal{I}_C^*)$, with $\mathcal{I}_B^* = \{\iota_1, \ldots, \iota_{(k+1)p}\}$ and $\mathcal{I}_C^* = \{\gamma_1, \ldots, \gamma_{(k+1)m}\}$, where $k$ denotes the maximum number of actuator-sensor-communication failures. In addition, let $\mathcal{I}_B^{(i)}$ and $\mathcal{I}_C^{(i)}$, with $i = 1, \ldots, k+1$, be the indices associated minimum collection of dedicated actuators and sensors $\bar{B}(\{\iota_i, \ldots, \iota_{ip}\})$ and $\bar{C}(\{\gamma_i, \ldots, \gamma_{im}\})$ respectively, where $(\bar{A}, \bar{B}(\{\iota_i, \ldots, \iota_{ip}\}))$ and $(\bar{A}, \bar{C}(\{\gamma_i, \ldots, \gamma_{im}\}))$ are structurally controllable and observable, respectively. Further, let $\mathcal{K}^{(i)}$ be the sparsest information pattern determined with Procedure 1 when the two collections $(\mathcal{I}_B^{(i)}, \mathcal{I}_C^{(i)})$ are considered. Then, under assumption **A1**, $(\mathcal{I}_B^*, \mathcal{I}_C^*, \mathcal{K}^*)$, where $\mathcal{K}^* = \bigcup_{j=1}^{k+1} \mathcal{K}^{(j)}$ is a solution to $\mathcal{P}_2$.

Furthermore, a solution to $\mathcal{P}_2$ can be determined in polynomial time, more precisely, $\mathcal{O}(k \max\{|\mathcal{V}| + |\mathcal{E}|, (|\mathcal{I}_B^*| + |\mathcal{I}_C^*|)^{\frac{5}{2}}\})$, where $|\mathcal{V}|$ and $|\mathcal{E}|$ are the number of vertices and edges of $\mathcal{D}(\bar{A}, \bar{B}(\mathcal{I}_B^*), \bar{C}(\mathcal{I}_C^*))$, respectively. ◇

### A. Resilient actuation-sensing-communication co-design for stabilization

In practice, it is often desirable to ensure only stabilizability of the closed-loop system, i.e., fixed modes can exist as long as they are stable. The same readily applies for the structurally (resilient) fixed modes. Specifically, we allow structurally fixed modes that are at the origin, and consequently we seek to placing the remaining ones, which we refer to as (possibly) *structurally unstable resilient fixed modes*. Specifically, consider the following variation of $\mathcal{P}_2$ which constraint is that of not having unstable structurally resilient fixed modes.

$\mathcal{P}_2^s$ Given $\bar{A}$, and a maximum of $k$ actuator-sensor-communication failures, we want to determine $(\bar{B}(\mathcal{I}_B^*), \bar{C}(\mathcal{I}_C^*), \bar{K}($ where $(\mathcal{I}_B^*, \mathcal{I}_C^*, \mathcal{K}^*)$ is a solution to the following optimization problem:

$$\begin{aligned} \min_{\mathcal{I}_B, \mathcal{I}_C, \mathcal{K}} \quad & |\mathcal{I}_B| + |\mathcal{I}_C| + |\mathcal{K}| \\ \text{s. t.} \quad & (\bar{A}, \bar{B}(\mathcal{I}_B \setminus \mathcal{J}_B), \bar{C}(\mathcal{I}_C \setminus \mathcal{J}_C); \bar{K}, \mathcal{J}) \text{ has no} \\ & \text{struct. unstable resilient fixed modes, for all} \\ & \mathcal{J} \subset \bar{\mathcal{K}}, \mathcal{J}_B \subset \mathcal{I}_B, \mathcal{J}_C \subset \mathcal{I}_C, \text{ and} \\ & |\mathcal{J}| + |\mathcal{J}_B| + |\mathcal{J}_C| \leq k, \end{aligned} \tag{5}$$

where $|\mathcal{I}_B|$ and $|\mathcal{I}_C|$ are the number of state variables actuated and measured, and $|\mathcal{K}|$ is the number of communication links established from the sensors to the actuators to obtain the decentralized closed-loop control system. ◦

In particular, it follows the feasibility conditions to $\mathcal{P}_2^s$ are attained if the following result holds.

*Corollary 1:* The structural system $(\bar{A}, \bar{B}, \bar{C}; \mathcal{K}, \{\Gamma_i\}_{i\in\Delta})$ has no structurally fixed modes (other than at the origin) if and only if in $\mathcal{D}(\bar{A}, \bar{B}, \bar{C}; \mathcal{K}, \Gamma_i) = (\mathcal{X} \cup \mathcal{U} \cup \mathcal{Y}, \mathcal{E}_{\mathcal{X},\mathcal{X}} \cup \mathcal{E}_{\mathcal{X},\mathcal{Y}} \cup \mathcal{E}_{\mathcal{U},\mathcal{X}} \cup (\mathcal{E}_{\mathcal{Y},\mathcal{U}} \setminus \mathcal{E}_{\Gamma_i}))$ for all $i \in \Delta$, with $\mathcal{E}_{\Gamma_i} = \{(y_i, u_j) : (j, i) \in \Gamma_i\}$, each state vertex $x \in \mathcal{X}$ is contained in an SCC which includes an edge of $\mathcal{E}_{\mathcal{Y},\mathcal{U}} \setminus \mathcal{E}_{\Gamma_i}$. ⋄

Thus, it is possible to waive the assumption stated in **A1**, while obtaining similar results to those in Theorem 4 that hold generally. Specifically, following the reasoning performed in the proof of Theorem 4, we obtain the following result.

*Corollary 2:* Let $\mathcal{I}_B$ and $\mathcal{I}_C$, be the indices associated minimum collection of dedicated actuators and sensors $\bar{B}(\mathcal{I}_B^*)$ and $\bar{C}(\mathcal{I}_C^*)$, with $\mathcal{I}_B^* = \{\iota_1, \ldots, \iota_{(k+1)p}\}$ and $\mathcal{I}_C^* = \{\gamma_1, \ldots, \gamma_{(k+1)m}\}$, and $k$ the maximum number of actuator-sensor-communication failures. In addition, let $\mathcal{I}_B^{(i)}$ and $\mathcal{I}_C^{(i)}$, with $i = 1, \ldots, k+1$, be the indices associated minimum collection of dedicated actuators and sensors $\bar{B}(\{\iota_i, \ldots, \iota_{ip}\})$ and $\bar{C}(\{\gamma_i, \ldots, \gamma_{im}\})$ respectively, where $(\bar{A}, \bar{B}(\{\iota_i, \ldots, \iota_{ip}\}))$ and $(\bar{A}, \bar{C}(\{\gamma_i, \ldots, \gamma_{im}\}))$ are structurally controllable and observable, respectively. Further, let $\mathcal{K}^{(i)}$ be the sparsest information pattern determined with Procedure 1 when the two collections $(\mathcal{I}_B^{(i)}, \mathcal{I}_C^{(i)})$ are considered. Then, under assumption **A1**, $(\mathcal{I}_B^*, \mathcal{I}_C^*, \mathcal{K}^*)$, where $\mathcal{K}^* = \bigcup_{j=1}^{k+1} \mathcal{K}^{(j)}$ is a solution to $\mathcal{P}_2^s$.

Furthermore, a solution to $\mathcal{P}_2^s$ can be determined in polynomial time, more precisely, $\mathcal{O}(k \max\{|\mathcal{V}| + |\mathcal{E}|, (|\mathcal{I}_B^*| + |\mathcal{I}_C^*|)^{\frac{5}{2}}\})$, where $\mathcal{V}$ and $\mathcal{E}$ are the number of vertices and edges of $\mathcal{D}(\bar{A}, \bar{B}(\mathcal{I}_B^*), \bar{C}(\mathcal{I}_C^*))$, respectively. ⋄

### *B. Stabilizing the closed-loop parameterized system*

Now, we discuss how the previous structural results can be leveraged to attain a specific gain matrix $K$ satisfying a specific information pattern, when the parameters of $(A, B, C)$ are considered, and such that

**ALGORITHM 1:** Stabilizability of the closed-loop system with possible communications between sensors and actuators given by $\mathcal{K}$, and a collection of possible communication failures $\{\Gamma_i\}_{i\in\Delta}$, with $\Gamma_i \subset \mathcal{K}$

**Input:** The system $(A, B, C)$, a set of possible communications between sensors and actuators $\mathcal{K}$ and a collection of possible communication failures $\{\Gamma_i\}_{i\in\Delta}$, with $\Gamma_i \subset \mathcal{K}$.

**Output:** a matrix $K$ such that $\tilde{A}(K)$ is Schur for $\{K \in [\bar{K}(\mathcal{K})]^{\Gamma_i}\}_{i\in\Delta}$.

1: Find feasible points $X_0, Y_0$ and $K_0$ that satisfy the constraints in (6). If a feasible point does not exist, then it is not possible to stabilize the system with such information pattern constraints;

2: At iteration $k > 0$, from $X_k$, $Y_k$ obtain the matrices $X_{k+1}$, $Y_{k+1}$ and $K_{k+1}$ solving

$$\begin{aligned} \min_{\{K\in[\bar{K}(\mathcal{K})]^{\Gamma_i}\}_{i\in\Delta},\ X,Y\in\mathbb{S}_+^N} \quad & \text{trace}(Y_k X_{k+1} + X_k Y_{k+1}) \\ \text{s.t.} \quad & \begin{bmatrix} X & \tilde{A}^\intercal(K) \\ \tilde{A}(K) & Y \end{bmatrix} \succ 0 \\ & \begin{bmatrix} X & I \\ I & Y \end{bmatrix} \succeq 0, \end{aligned}$$

where $\tilde{A}(K) = \begin{bmatrix} A & B \\ C & K \end{bmatrix}$.

3: If the matrix

$$\tilde{A}_{k+1} = \begin{bmatrix} A & B \\ C & K_{k+1} \end{bmatrix}$$

is Schur for $K_{k+1} \in [\bar{K}]^{\Gamma_i}, i \in \Delta$, stop the algorithm. Otherwise, set $k = k + 1$ and go to Step 2.

the closed-loop system is always asymptotically stable even if some failure occurs. Consequently, we are left to ensure that $K$ is such that the eigenvalues of $A + BKC$ lie within the unit circle in the complex plane for $K \in [\bar{K}]^{\Gamma_i}_{p,m}$ for $i \in \Delta$. This condition is equivalent to say that $\tilde{A}(K) = \begin{bmatrix} A & B \\ C & K \end{bmatrix}$ is Schur. Subsequently, we obtain the following result.

*Theorem 5:* Given $(A, B, C)$, a set of possible communications between sensors and actuators $\mathcal{K}$ and a collection of possible communication failures $\{\Gamma_i\}_{i\in\Delta}$, with $\Gamma_i \subset \mathcal{K}$, the $N \times N$ matrix $\tilde{A}(K)$ is Schur for $\{K \in [\bar{K}(\mathcal{K})]^{\Gamma_i}\}_{i\in\Delta}$, if and only if the following optimization problem

$$\begin{aligned} \min_{\{K\in[\bar{K}(\mathcal{K})]^{\Gamma_i}\}_{i\in\Delta},\ X,Y\in\mathbb{S}_+^N} \quad & \text{trace}(XY) \\ \text{s.t.} \quad & \begin{bmatrix} X & \tilde{A}^\intercal(K) \\ \tilde{A}(K) & Y \end{bmatrix} \succ 0, X = Y^{-1} \end{aligned} \tag{6}$$

is feasible with optimal cost $N$.

◇

*Remark 5:* Theorem 5 characterizes the existence of unstable resilient fixed modes versus the non-existence of resilient structural fixed modes. More specifically, if there is a $K$ satisfying (6) than it follows that the system has no unstable resilient fixed modes. Conversely, if the system has unstable resilient fixed modes, then there is no feasible solution to (6). $\diamond$

Theorem 5 provides a characterization of the solutions, but it is non convex. Hence, we propose to approximate it by a recursive scheme as described in Algorithm 1, where each iteration can be efficiently solved using standard LMI toolboxes, and whose feasibility properties are described in the following result.

*Theorem 6:* Algorithm 1 determines a Schur matrix $\tilde{A}(K)$ for $\{K \in [\bar{K}(\mathcal{K})]^{\Gamma_i}\}_{i\in\Delta}$, if the sequence trace$(Y_k X_{k+1} + X_k Y_{k+1})$ converges to $2N$. $\diamond$

## V. Illustrative example

In this section, we illustrate the main results obtained in the context of the power electric grid. More specifically, we consider the IEEE 5-bus system as modeled in [39] that consists in a linearized model under normal operating conditions, and provide the solution to problem $\mathcal{P}_2$ when $k = 1$ actuator-sensor-communication failure can occur. In addition, we determine the gain (for a specific realization of the system) subject to information pattern constraints as proposed in Algorithm 1.

The IEEE 5-bus power system consists of 3 synchronous generators connected to 2 aggregate loads coupled through the network topology. The state digraph $\mathcal{D}(\bar{A})$ associated with its linearized dynamics $\bar{A}$ is depicted in Figure 2, and it consists of 18 state variables, which physical interpretation is as follows: $x_2, x_5$ and $x_8$ represents the mechanical power of the turbine of the generator 1-3, respectively. The generators' frequency is captured by $x_1, x_4$ and $x_7$ for generator 1-3 respectively, and $x_3, x_6$ and $x_9$ their valve opening. In addition, $x_{11}$ and $x_{13}$ is the real energy consumed by the aggregate load 1 and 2 respectively, and $x_{10}$ and $x_{12}$ their frequency measured locally. The different components are connected through the injected/received power to/from the network at the connection site, which dynamics depend on the frequency of the components on the neighboring buses; the injected/received power variables for generator 1-3 and load 1-2 are captured by $x_{14}, x_{15}, x_{16}, x_{17}$ and $x_{18}$, respectively. Furthermore, notice that in the present setup, the assumption of considering dedicated inputs/outputs is consistent with the nature of the physical quantities the state variables aim to capture.

We start by noticing that the state digraph depicted in Figure 2 consists of three SCCs, where one is non-bottom linked and the other two are non-top linked SCCs. Let us denote the not-top linked SCCs containing the state variable $x_{11}$ (belonging to load 1) and $x_{13}$ (belonging to load 2) by $\mathcal{N}_1^T$ and $\mathcal{N}_2^T$, respectively. In addition, let us denote the non-bottom linked SCC by $\mathcal{N}_1^B$. Furthermore, notice that assumption **A1** holds, since each state variable (vertex) has a self-loop, hence, the collection of these edges correspond to a perfect matching of the state bipartite graph. Therefore, from Lemma 2 and its dual, we have that $(\mathcal{I}_B^{(1)} \equiv \{1,2\}, \mathcal{I}_C^{(1)} \equiv \{1\})$ and $(\mathcal{I}_B^{(2)} \equiv \{3,4\}, \mathcal{I}_C^{(2)} \equiv \{2\})$ correspond to two possible minimum collections of dedicated inputs and outputs ensuring structural controllability and observability, respectively. In addition, $\mathcal{I}_B^* = \{11, 13, 11, 13\}$ and $\mathcal{I}_C^* = \{10, 12\}$, which implies that $\bar{B}(\mathcal{I}_B^*) = [e_{11}\ e_{13}\ e_{11}\ e_{13}]$ and $\bar{C}(\mathcal{I}_C^*) = [e_{10}\ e_{12}]^\intercal$, where $e_i \in \mathbb{R}^{18}$ denotes the $i^{\text{th}}$ canonical vector with 1 in entry $i$ and zero elsewhere, ensure structural controllability and observability when one dedicated actuator-sensor fails, recall Lemma 2 and its dual. Notice that from a physical point of view the actuation of the variables $x_{11}$ and $x_{13}$ means that we just need to be able to actuate the real power consumed by the aggregate load. In particular, the actuation capabilities considered assume that the generators inject electric power at a constant rate, delegating the responsibility of stabilizing the network to the loads; in other words, the loads can both receive and inject electric power in the network, playing the role of *prosumers*. It is also worth noting that the state variables considered to be measured were the frequency of the aggregate load since it is locally inferred by one of the loads, hence, requiring effective communication over a network of one data package to the remaining load. Nonetheless, any other state variable in the non-bottom linked SCC would be a viable option.

Now, consider Procedure 1 when $(\mathcal{I}_B^{(1)}, \mathcal{I}_C^{(1)})$ is considered. Then, $\mathcal{B} \equiv \mathcal{B}(\mathcal{I}_B^{(1)}, \mathcal{I}_C^{(1)}, \mathcal{E}_{\mathcal{I}_B^{(1)}, \mathcal{I}_C^{(1)}})$ where $\mathcal{E}_{\mathcal{I}_B^{(1)}, \mathcal{I}_C^{(1)}} = \{(1,1), (2,1)\}$, i.e., the the dedicated actuators controlling the state variables in the non-top linked SCC $\mathcal{N}_1^T$ and $\mathcal{N}_2^T$ reach the dedicated sensor measuring the state variable in the non-bottom linked SCC $\mathcal{N}_1^B$. A possible maximum matching $M^*$ associated with $\mathcal{B}$ is $M^* = \{(2,1)\}$, which implies that $\mathcal{K} = \{(1,2)\}$ is a possible sequential-pairing. In particular, notice that $M^* \cup \mathcal{K} = \{(1,2), (2,1)\}$ forms a cycle. Subsequently, $\mathcal{D}(\bar{A}, \bar{B}(\mathcal{I}_B^*), \bar{C}(\mathcal{I}_C^*), \bar{K}(\mathcal{K}))$ contains only one non-top and one non-bottom linked SCC containing state variables. Therefore, we just need to consider $\mathcal{K}' = \{(1,1)\}$ as prescribed in Procedure 1, and it follows that $\mathcal{D}(\bar{A}, \bar{B}(\mathcal{I}_B^*),\ \bar{C}(\mathcal{I}_C^*), \bar{K}(\mathcal{K}^{(1)} \equiv \mathcal{K} \cup \mathcal{K}'))$ has a single SCC containing state variables with two feedback links on it.

Similarly, Procedure 1 can be used when the pair $(\mathcal{I}_B^{(2)}, \mathcal{I}_C^{(2)})$ is considered, which leads to $\mathcal{K}^{(2)} = \{(2,3), (2,4)\}$. The digraph $\mathcal{D}(\bar{A}, \bar{B}(\mathcal{I}_B^*), \bar{C}(\mathcal{I}_C^*), \bar{K}(\mathcal{K}^* \equiv \mathcal{K}^{(1)} \cup \mathcal{K}^{(2)}))$ that accounts for the obtained solution is depicted in Figure 2.

In addition, by noticing that all state variables in $\mathcal{D}(\bar{A})$ have self-loops, it follows that assumption **A1** holds. Hence, invoking Theorem 4, we obtain that $(\mathcal{I}_B^*, \mathcal{I}_C^*, \mathcal{K}^*)$ is a solution to $\mathcal{P}_2$ when $k = 1$.

In addition, we notice that by Corollary 2 the same procedure ensures the existence of no structurally fixed mode other than those at the origin, even if assumption **A1** does not hold.

Now, consider Algorithm 1 when all generators and loads have equal parameterization, perturbed by zero mean Gaussian noise with standard deviation squared equal to $0.1$. In particular, the nonzero entries of $A$ are as follows: $A_{1,1} = 0.2309$, $A_{3,1} = 0.8005$, $A_{14,1} = 0.7662$, $A_{15,1} = 0.9755$, $A_{17,1} = 1.0123$,

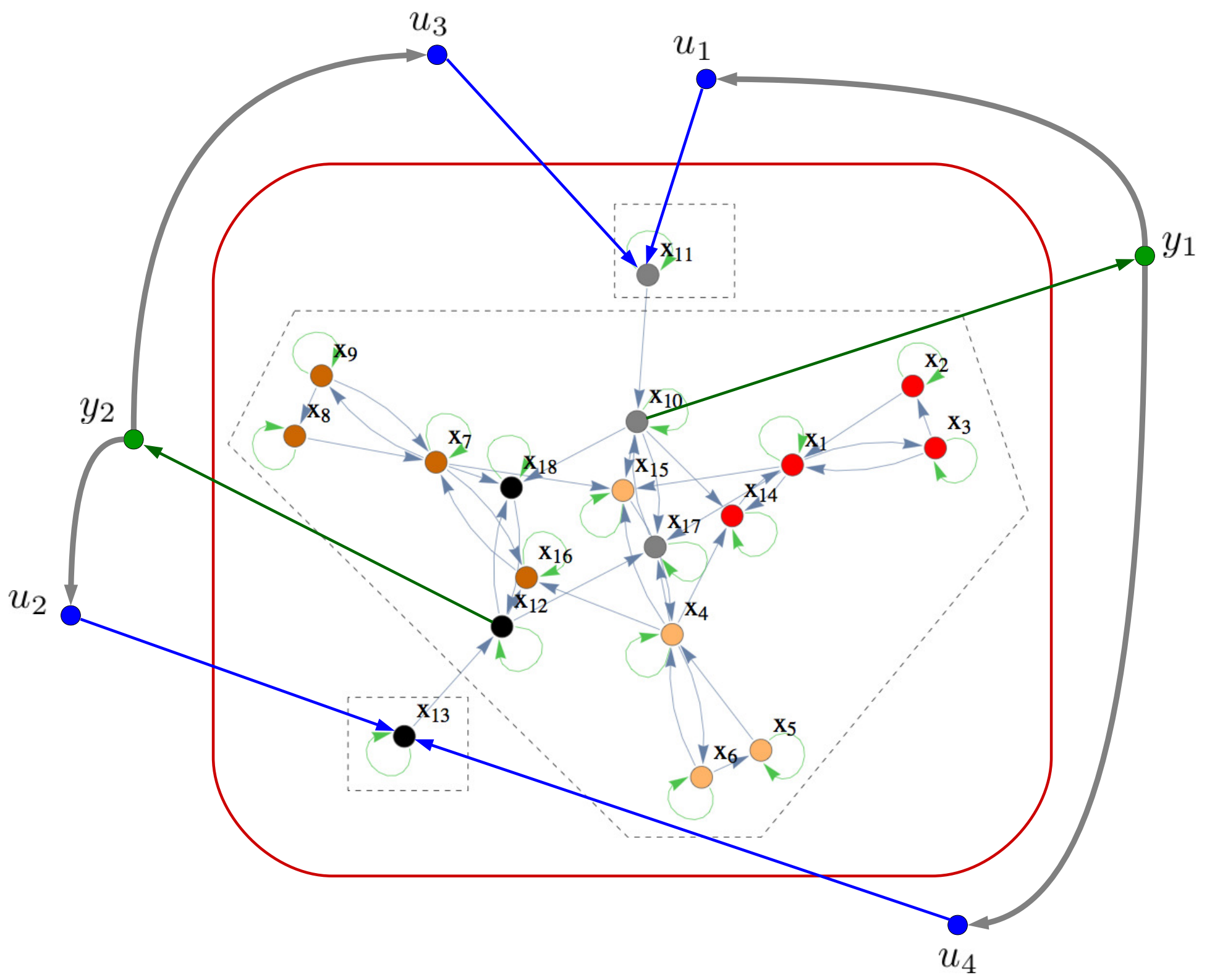

Fig. 2: This figure presents the digraph $\mathcal{D}(\bar{A},\bar{B},\bar{C},\bar{K}(\mathcal{K}^*))$, with four dedicated actuators depicted in blue, two dedicated sensors depicted in dark green, and the feedback links depicted by gray arrows. The state digraph representation of the 5-bus system is enclosed by the red round box, where the state vertices in gray and black belong to load 1 and 2, respectively. The state vertices in red, orange and brown belong to generator 1, 2 and 3, respectively. The different SCCs are enclosed the dashed boxes, where the squared dashed boxes represent the non-top linked SCCs and the remaining SCC is non-bottom linked. Finally, notice that all state variables have self-loops depicted by green edges. The physical interpretation of the state variables is provided in the main text.

$A_{1,2}=0.6275$, $A_{2,2}=0.1166$, $A_{1,3}=0.4365$, $A_{2,3}=-0.1303$, $A_{3,3}=-0.0100$, $A_{4,4}=0.2154$, $A_{6,4}=0.8283$, $A_{14,4}=1.1273$, $A_{15,4}=0.7095$, $A_{16,4}=1.0272$, $A_{17,4}=1.1220$, $A_{4,5}=0.8733$, $A_{5,5}=0.1938$, $A_{4,6}=0.3790$, $A_{5,6}=-0.0800$, $A_{6,6}=0.2679$, $A_{7,7}=0.1458$, $A_{9,7}=1.1761$, $A_{15,7}=0.7831$, $A_{16,7}=0.6320$, $A_{18,7}=0.9408$, $A_{7,8}=0.6544$, $A_{8,8}=0.1681$, $A_{7,9}=0.3592$, $A_{8,9}=-0.2192$, $A_{9,9}=0.0637$, $A_{10,10}=0.4795$, $A_{14,10}=0.8029$, $A_{15,10}=0.8851$, $A_{17,10}=0.5781$, $A_{18,10}=1.0721$, $A_{10,11}=-0.0275$, $A_{11,11}=0.7626$, $A_{12,12}=0.2686$, $A_{16,12}=0.7851$, $A_{17,12}=1.0285$, $A_{18,12}=0.4553$, $A_{12,13}=0.1462$, $A_{13,13}=0.7775$, $A_{1,14}=0.0491$, $A_{14,14}=0.4760$, $A_{4,15}=0.2080$, $A_{15,15}=0.5730$, $A_{7,16}=-0.0535$, $A_{16,16}=0.5155$, $A_{10,17}=0.0492$, $A_{17,17}=0.0407$, $A_{12,18}=-0.0668$, and $A_{18,18}=0.2072$. The eigenvalues of $A$ have norms given as follows: 1.0072, 0.7775, 0.7626, 0.6652, 0.6636, 0.6636, 0.6017, 0.5674, 0.5674, 0.4802, 0.4802, 0.4678, 0.4081, 0.3596, 0.3001, 0.3001, 0.1359, and 0.1359. Hence, $A$ is not stable since one eigenvalue lie outside the unitary circle in the complex plane. In addition, consider the set of communication links $\mathcal{K}=\mathcal{K}^{(1)}\cup\mathcal{K}^{(2)}$, and assume that any communication link can fail or be compromised by an attacker, more precisely, $\Gamma_1=\{(1,1)\}$, $\Gamma_2=\{(1,2)\}$, $\Gamma_3=\{(2,3)\}$ and

$\Gamma_4 = \{(2,4)\}$. Consequently, by performing Algorithm 1, we obtain the following gain:

$$K = \begin{bmatrix} 0.6445 & 0 \\ -0.3043 & 0 \\ 0 & -1.0079 \\ 0 & -0.0000043329 \end{bmatrix},$$

where largest norm associated with the eigenvalue of the closed-loop system is $0.9918$, which implies that all eigenvalues lie within the unit circle in the complex plane; hence, the closed-loop system is stable. In addition, if an actuator-sensor-communication failure encoded by $\Gamma_i$ $(i = 1, \ldots, 4)$ occurs, then the gains used by the closed-loop system are effectively modified due to loss of specific sensor-actuator feedback links, i.e., due to the failure, the closed-loop system effectively uses the following gains:

$$K^{\Gamma_1} = \begin{bmatrix} 0 & 0 \\ -0.3043 & 0 \\ 0 & -1.0079 \\ 0 & -0.0000043329 \end{bmatrix}, \quad K^{\Gamma_2} = \begin{bmatrix} 0.6445 & 0 \\ 0 & 0 \\ 0 & -1.0079 \\ 0 & -0.0000043329 \end{bmatrix},$$

$$K^{\Gamma_3} = \begin{bmatrix} 0.6445 & 0 \\ -0.3043 & 0 \\ 0 & 0 \\ 0 & -0.0000043329 \end{bmatrix}, \quad K^{\Gamma_4} = \begin{bmatrix} 0.6445 & 0 \\ -0.3043 & 0 \\ 0 & -1.0079 \\ 0 & 0 \end{bmatrix},$$

where the norm of the largest eigenvalue of the closed-loop system are $0.9981$, $0.9983$, $0.9989$, and $0.9918$, respectively. Thus, the closed-loop system remains stable under any of these effective gains, i.e., under any of the considered set of attack conditions.

*Computational Implementation Details:* We notice that the gain matrices have to be square due to the LMI constraint $\begin{bmatrix} X & I \\ I & Y \end{bmatrix} \succeq 0$, so we considered two additional sensors that do not measure any state variable, i.e., they correspond to columns of $\bar{C}$ with only zeros. For simulation purposes, we initialized $X_0$ and $Y_0$ to be $22 \times 22$ identity matrices. Therefore, as prescribed in Theorem 6 the gains to stabilize the closed-loop system exist if the cost function of the iterative algorithm (Algorithm 1) converges to $44$. We resort to the CVX (as LMI solver) and Matlab to execute Algortihm 1, which terminated after approximately 50 minutes in an Apple laptop with the following specifications: 2.6 GHz Intel Core i5 with 8 GB (1600 MHz DDR3) of memory, yet after approximately 25 min the objective of attaining a Schur matrix was already attained. $\circ$

## VI. Conclusions and further research

In this paper, we provided a solution to the problem of determining the minimum number of actuators, sensors and communication between these, for linear discrete time-invariant plants, and such that desirable control objectives of the closed-loop system are achieved. Namely, those achieving arbitrary performance in terms of the system response and/or stabilizability with static output feedback, under possible disruptive scenarios such as natural failures, or attacks by malicious agents, of actuators, sensors and communication links.

More precisely, we introduced the notion of resilient fixed-modes free system that ensures the non-existence of fixed modes when the actuation-sensing-communication structure is compromised; hence, crucial to ensure the proper functioning of the system under failures and/or attacks. In addition, we provided a graph-theoretical representation that ensures almost always the non-existence of resilient fixed modes, which we used to address the problem of designing the minimum actuation-sensing-communication structure that ensures a system without resilient fixed modes almost always. Furthermore, we discussed how the assumptions of the dynamical system structure can be waived in the context of stabilization, while the methodology presented holds. Finally, we provided a convex optimization algorithm that ensures stabilizability of the closed-loop system under possible actuator-sensor-communication failures, and we illustrated the results in the context of smart grid, using the IEEE 5-bus system.

Future research will address how the proposed methodology can be extended to handle intermittent communication schemes and delays of data, as well as considering potential actuation-sensing-communication costs. In addition, we will focus on developing more efficient and general methods to find the gains satisfying information patterns ensuring the closed-loop system have no resilient fixed modes.

## APPENDIX

Proof of Theorem 1: The proof follows by comparison with the definition of structurally fixed modes and its characterization provided in Theorem 4 in [44]. More precisely, from the definition of resilient fixed modes (Definition 1) it follows that $\Sigma(A, B, C; \mathcal{K}, \{\Gamma_i\}_{i\in\Delta}) = \emptyset$ if and only if $\bigcap_{K\in[\bar{K}(\mathcal{K})]^{\Gamma_i}_{p,m}} \sigma(A+BKC) = \emptyset$ for all $i \in \Delta$. Subsequently, by definition of structurally fixed mode, it follows that a system has no structurally resilient fixed modes with respect to an information pattern $\bar{K}(\mathcal{K})$ and a set of possible failures $\{\Gamma_i\}_{i\in\Delta}$ if and only if $(\bar{A}, \bar{B}, \bar{C}; \mathcal{K}, \Gamma_i)$ has no structurally fixed modes for all $i \in \Delta$. Thus, invoking Theorem 4 in [44] the result follows. ■

Proof of Theorem 2: The proof follows by reducing the minimum $k$ robust structural controllability problem, i.e., the problem of determining the minimum number of dedicated inputs required to ensure structural controllabiliy under at most $k > 1$ input failures and given the structrural dynamical matrix $\bar{A}$, to the resilient co-design problem $\mathcal{P}_2$. Specifically, the same parameters $\bar{A}$ and $k > 1$ can be adopted for the problem $\mathcal{P}_2$, and a solution to $\mathcal{P}_2$ provides in particular, a solution to the minimum $k > 1$ robust structural controllability problem can be retrieved from a solution to $\mathcal{P}_2$ since the system has no structurally fixed modes only if it is structurally controllable and structurally observable [27]. Subsequently, if $k$ failures occur in the inputs, it follows that $\mathcal{I}_B$ must contain a solution to the minimum $k$ robust structural controllability problem. In fact, we notice that $\mathcal{I}_B$ consists into a solution to the latter problem, due to the equivalence between actuation-sensing-communication failure in our modeling setup see Figure 1. Hence, the result follows, i.e., $\mathcal{P}_2$ is at least as difficult as solving the minimum $k$ robust structural controllability problem, which is known to be NP-hard [34]. ■

Proof of Lemma 1: First, notice that by construction Procedure 1 provides a feasible solution. To show optimality, we proceed as follows: (i) all dedicated actuators and sensors in the minimum collections ensuring structural controllability and observability are required; and (ii) we show that we need at least $|\mathcal{I}_C^*|$ (assuming $|\mathcal{I}_C^*| \geq |\mathcal{I}_B^*|$) feedback links. First, observe that (i) follows from noticing that otherwise there would exist a non-top or non-bottom linked (in case structural controllability/observability is compromised, respectively), hence, condition (a) in Theorem 1 does not hold. Secondly, all dedicated actuators and sensors must be used to obtain a feasible information pattern, since if there is one dedicated actuator (respectively, dedicated sensor) that is not used then, the non-top linked SCC (respectively, non-bottom linked SCC) that contains the state variable it controls (respectively, measures) will stand alone when the digraph associated with the closed-loop system is considered; hence, this SCC has no feedback link on it, which violates condition (a) in Theorem 4 in [44]. In summary, we need to use $|\mathcal{I}_C^*|$ feedback links, i.e., $|\mathcal{K}^*| \geq |\mathcal{I}_C^*|$, and the result follows.

Finally, the computational complexity is due to the following components: $\mathcal{E}_{\mathcal{I}_B,\mathcal{I}_C}$ can be determined through a depth-first search rooted in the dedicated input $i \in \mathcal{I}_B$, which can be achieved in $\mathcal{O}(|\mathcal{V}|+|\mathcal{E}|)$. In addition, the sequential pairing requires the computation of the maximum matching associated with $\mathcal{B}$ that can be computed using the Hungarian algorithms, which runs in $\mathcal{O}(\sqrt{n}m)$ for a bipartite with $n$ vertices and $m$ edges, so we obtain $(|\mathcal{I}_B^*|+|\mathcal{I}_C^*|)^{\frac{5}{2}}$. Thus, because all the remaining operations have linear-time complexity, the result follows. ■

Proof of Theorem 3: The feasibility follows from noticing that each dedicated actuator and sensor must be involved in at least $k+1$ feedback links. More precisely, suppose by contradiction that this is not the case, then if a dedicated input has $k$ (or less) feedback links, then when these $k$ feedback links fail the dedicated input does not have any feedback link incoming on it. Therefore, the non-top linked SCC that contains the state variable it actuates does not belong to an SCC with a feedback link on it, which violates condition (a) in Theorem 4 in [44]. Similar reasoning applies when dedicated sensors are considered. Thus, each dedicated actuator and sensor must be involved in at least $k+1$ feedback links, which is ensured by considering $\mathcal{K}^* = \mathcal{K}^{(1)} \cup \ldots \cup \mathcal{K}^{(k+1)}$, where $\mathcal{K}^{(i)}$ is determined using Procedure 1 and $\mathcal{K}^{(i)} \cap \mathcal{K}^{(j)} = \emptyset$ for $i \neq j$. In addition, the optimality, i.e., the fact that $\bar{K}(\mathcal{K}^*)$ is a sparsest information pattern ensuring that the closed-loop system has no structurally resilient fixed modes, follows from noticing that each dedicated actuator and sensor is involved in exactly in $k+1$ feedback links.

The computational complexity follows the same steps as in Lemma 1, by noticing that the procedure is executed $k$ times, having in consideration Remark 4. ■

Proof of Lemma 2: The proof follows by invoking Theorem 3 in [27], and noticing that there must exists at least one state variable actuated in each non-top linked SCC after the failures occur. ■

Proof of Theorem 4: The feasibility is obtained as follows: (i) Lemma 2 and its dual, prescribe the existence of $k{+}1$ minimum collections of dedicated actuators and sensors ensuring structural controllability and observability, respectively; (ii) for each minimum collection $i$ of dedicated actuators and sensors, Theorem 3 ensures that $\bar{K}(\mathcal{K}^{(i)})$ is the sparsest information pattern ensuring a closed-loop system without structurally fixed modes; and (iii) because there exist $k+1$ digraphs verifying Theorem 4 in [44], i.e., $\mathcal{D}(\bar{A},\bar{B},\bar{C},\bar{K}(\mathcal{K}^{(i)}))$ for $i=1,\ldots,k+1$, it follows that Theorem 1 holds. Furthermore, the optimality follows by noticing that minimum of the sum of the above quantities is achieved by the proposed design.

The computational complexity follows by reproducing the same steps as in Theorem 3 and noticing that determining $k$ minimum collection of dedicated actuators and sensors, which requires only to identify the non-top and non-bottom linked SCCs, which can be achieved by determining the decomposition of the digraph $\mathcal{D}(\bar{A})$ into SCCs that can be implemented in $\mathcal{O}(|\mathcal{V}|+|\mathcal{E}|)$. ■

Proof of Theorem 5: From Theorem 2-(i) in [45], we have that a matrix $\tilde{A}(K)$ is Schur if and only if there exist symmetric matrices $X$ and $Y$ such that $\begin{bmatrix} X & \tilde{A}(K)^{\intercal} \\ \tilde{A}(K) & Y \end{bmatrix} \succ 0$, and $X=Y^{-1}$. In addition, $X=Y^{-1}$ if and only if $X,Y$ are a solution to the following optimization problem

$$\begin{aligned} \min_{X,Y\in\mathbb{S}^n_+} \quad & \text{trace}(XY) \\ \text{s.t.} \quad & \begin{bmatrix} X & I \\ I & Y \end{bmatrix} \succeq 0 \end{aligned}$$

and the optimal cost is $n$ [46], and the result follows. ■

Proof of Theorem 6: The results follows by invoking the algorithm proposed in [46]. In particular, the iterative method with cost trace$(Y_kX_{k+1}{+}X_kY_{k+1})$ will always converge, such that the constraint $Y=X^{-1}$ satisfied if the cost converges to $2N$ [47].

■